\documentclass[11pt,a4paper,Color]{article}
\usepackage{amssymb, paralist}
\usepackage{amsfonts, amsmath, wasysym}
\usepackage{pdfsync}
\usepackage{latexsym}
\usepackage{graphicx, color}
\usepackage{indentfirst}
\usepackage{bm}

\usepackage[draft]{optional}

\usepackage{fancyhdr}
\newcommand{\mc}[1]{\mathcal{#1}}

\input xypic.sty

\newcommand{\com}[1]{\opt{draft}{\textcolor{red}{
$\LHD$ #1 $\RHD$\marginpar{\textcolor{red}{$\blacksquare$}}}}}

\def\qed{\hfill {\large ${\sqcup\!\!\!\!\sqcap}$}}

\newenvironment{demo}{{\bf Proof }}
{\qed \\}

\newcommand{\ba}{\begin{array}}
\newcommand{\ea}{\end{array}}

\newcommand{\re}{\mathbb R}

\newcommand{\flecha}{\longrightarrow}
\newcommand{\<}{\left<}
\renewcommand{\(}{\left(}
\newcommand{\lb}{\label}
\newcommand{\nn}{\nonumber}
\newcommand{\fracc}{\displaystyle\frac}

\renewcommand{\>}{\right>}
\renewcommand{\)}{\right)}
\newcommand{\eps}{\ensuremath{\varepsilon}}
\newcommand{\gor}{\ensuremath{\widetilde}}

\newcommand{\us}{\underset}
\def\bal{\begin{align}}
\def\eal{\end{align}}

\def\z{{\frak z}}
\def\de{{\frak d}}

\newcommand{\bde}{\begin{defi}}
\newcommand{\ede}{\end{defi}}

\numberwithin{equation}{section}
\def\be{\begin{equation}}
\def\ee{\end{equation}}

\def\og{{\overline g}}
\def\oH{{\overline H}}
\def\oR{{\overline R}}
\def\oRic{{\overline Ric}}
\def\oM{{\overline M}}
\def\oN{{\overline \nabla}}

\def\wg{{\widetilde{\frak g}}}
\def\fg{{\frak g}}

\def\tM{{\widehat M }}

\def\tg{{\widehat g }}

\def\Sn{{\mathbb{S}^{n-1}}}

\def\a{\alpha}
\def\p{\varphi}

\def\a{\alpha}

\def\ral{\sqrt{|\lambda|}}
\def\vtr{{\rm tr}}

\def\nablaa{\overline{\nabla}}

\def\vle{{\rm vol}}
\def\parcial#1#2{\frac{\partial #1}{\partial#2}}

\def\deri#1#2{\frac{d #1}{d#2}}

\def\flecha{\longrightarrow}
\def\fle{\rightarrow}

\def\vt{\frak{t}}
\def\ds{\displaystyle}

\newtheorem{defi}{\hspace{12pt} Definition}
\newtheorem{teor}{\hspace{12pt} Theorem}
\newtheorem{prop}[teor]{\hspace{12pt} Proposition}
\newtheorem{notacion}[defi]{\hspace{12pt} Notation}
\newtheorem{lema}[teor]{\hspace{12pt} Lemma}

\newtheorem{nota}{\hspace{12pt} Remark}
\newtheorem{coro}[teor]{\hspace{12pt} Corollary}

\numberwithin{lemap}{teor}
\numberwithin{corop}{teor}

\numberwithin{ejer}{subsection}
\numberwithin{ejemplo}{subsection}

\begin{document}

\pagestyle{myheadings}\markboth{E. Cabezas-Rivas and V. Miquel}{Volume-preserving MCF between two equidistants}



\title{\vspace{-2cm} Volume preserving mean curvature flow of revolution hypersurfaces between two equidistants }

\author{ Esther Cabezas-Rivas \medskip \\
\footnotesize{ Westf.~Wilhelms Universit\"at M\"unster, Mathematisches Institut,} \\
     \footnotesize{ Einsteinstr.~62, 48149 M\"unster (Germany)}\\
     \footnotesize{Tel.~+49 251 83-33748, Fax +49 251 83-32711}\\
     \footnotesize{email: E.Cabezas-Rivas@uni-muenster.de}\bigskip\\
Vicente Miquel \medskip \\
   \footnotesize{Universidad de Valencia }\\
    \footnotesize{ 46100-Burjassot (Valencia) Spain} \\
    \footnotesize{email: miquel@uv.es }}


\maketitle

\vspace{-1cm}
\begin{abstract}
In a rotationally symmetric space $\oM$ around an axis $\mc A$ (whose precise definition includes all real space forms), we consider a domain $G$ limited by two equidistant hypersurfaces orthogonal to $\mc A$. Let $M \subset \oM$ be a  revolution hypersurface  generated by a graph over $\mc A$, with boundary in $\partial G$ and orthogonal to it. We study the evolution $M_t$ of $M$ under the volume-preserving mean curvature flow requiring that the boundary of $M_t$  rests on $\partial G$ and keeps orthogonal to it.  We prove that: a)  the generating curve of $M_t$ remains a graph; b) the flow 
exists while $M_t$ does not touch the axis of rotation;  c)  under a suitable hypothesis relating the enclosed volume and the area  of $M$, the flow is defined for every  $t\in [0,\infty[$ and a sequence of hypersurfaces $M_{t_n}$ converges to a   revolution hypersurface of constant mean curvature. Some key points are: i) the results are true even for ambient spaces with positive curvature, ii) the averaged mean curvature does not need to be positive and iii) for the proof it is necessary to carry out a detailed study of the boundary conditions.
\end{abstract}

{\bf Mathematics Subject Classification (2010)} 53C44

\section{Introduction and Main Results }\lb{In}

\subsection{Background about volume preserving evolution}
 A family of immersions $X_t: M \flecha \oM$, $t\in[0,T[$, of a $n$-dimensional compact manifold $M$ into a $(n+1)$-dimensional Riemannian manifold $(\oM,\og)$ is called a {\it Volume Preserving Mean Curvature Flow} ({\sc vpmcf})  if it is  a solution of the equation
\begin{equation}\label{vpmf}
\parcial{X_t}{t} = (\oH_t- H_t)\ N_t,
\end{equation}
where $\oH_t$ is the averaged mean curvature
$
\oH_t =\ds\frac{ \int_{M} H_t d \mu_t }{|M_t| }
$
of the immersion $X_t$, being $ d \mu_t$  the canonical volume element of the Riemannian manifold  $M_t=(M,X_t^*\og)$, $|M_t|$ its $n$-volume (which we shall call \lq\lq area''), $N_t$ the unit normal vector field pointing outward (if each $X_t(M)$ encloses a domain $\Omega_t$) and $H_t$ the mean curvature of the immersion $X_t$, with the following conventions: the Weingarten map $L_t:TM\flecha TM$ is given by $L_t Z = \oN_Z N_t$ and $H_t$ is its trace. Sometimes we shall also use the notation  $M_t=X_t(M)$. 

The presence of the global term $\oH$ in  equation \eqref{vpmf} has two major consequences: {\bf a)} keeps the enclosed volume constant while the area decreases, and {\bf b)} makes the usual techniques in geometric flows (e.g. the application of maximum principles) either fail or become more subtle. The resultant evolution problem is particularly {\it appealing} -since from a) it is specially well suited for applications to the isoperimetric problem- and {\it challenging} because b) causes a a plethora of extra complications; for instance, a basic principle for the ordinary mean curvature flow ({\it the comparison principle}) fails in general for \eqref{vpmf}, e.g., an initially embedded curve may develop self-intersections (cf. \cite{MaSi}). Hence the present knowledge of this flow is considerably poorer than that of the unconstrained evolution.

The {\sc vpmcf} has been studied under convexity assumptions for an initial closed hypersurface either within a Euclidean or Hyperbolic ambient space (cf. \cite{Hu87} and \cite{CaMi1}, respectively). There is intuitive evidence, as pointed out by G. Huisken in \cite{Hu87}, that the preservation of convexity may fail in ambient manifolds with positive curvature. One can also find stability results: if the initial hypersurface is close enough to a model constant mean curvature hypersurface, then it flows to one model (see \cite{AF}, \cite{ES}, \cite{CaMi1} and \cite{Li}).

After dealing with convexity assumptions, it is natural to wonder whether there is another natural geometric condition, invariant under \eqref{vpmf}, which still softens the problems caused by the global term. A good choice seems to take the initial $M$ to be a revolution hypersurface generated by the graph of a function over the axis of revolution of $M$. This was done for the Euclidean space in  \cite{Ath1, Ath2}. Later on, in \cite{CaMi2}, we considered $M$ within a wider family of ambient spaces (including the Euclidean and the Hyperbolic ones) for which still makes sense the notion of revolution hypersurface around an axis $\mathcal A$. 

The papers \cite{Ath1, Ath2, CaMi2} study the evolution under \eqref{vpmf} of $M$ as above, with boundary intersecting  orthogonally two totally geodesic hypersurfaces $\pi_{tg}$ orthogonal to $\mathcal A$, and  requiring that the evolving hypersurface meets $\pi_{tg}$ orthogonally at each time. When $\oM$ is not Euclidean, negativity of some of its sectional curvatures is imposed. It is proved:

{\it
\begin{quote}
\begin{itemize}
\item[ A] While the evolving hypersurface does not touch $\mathcal{A}$, the flow exists and  the generating curve remains a graph over $\mathcal{A}$.

\item[ B] Under a hypothesis relating the enclosed volume to the area of $M$, we achieve  long  time existence, and
convergence of a sequence $M_{t_n}$ to a revolution hypersurface of constant mean curvature.   
\end{itemize}
\end{quote}
}

In the Euclidean space, the hypersurfaces $\pi_{tg}$ are parallel hyperplanes, so they are at constant distance from each other; however, this is not any more true in the more general ambient spaces studied in \cite{CaMi2}. Then it is natural to address the same problem, but considering regions limited by hypersurfaces at constant distance.

 The main concern of the present paper will be the proof of the statements corresponding to A and B when changing $\pi_{tg}$ by equidistant limiting hypersurfaces. To understand some interesting issues arising in the new setting (cf. Section 1.3), it is important to highlight the following facts about the proofs of A and B in  \cite{Ath1,Ath2,CaMi2}.
\begin{itemize}
\item[(1)]  An isometry of the ambient space allows  to extend the problem to another bigger domain with symmetry so that the original boundary points become interior points and the maximum principle applies. Accordingly, the boundary of the evolving hypersurface does not cause any extra complication.

\item[(2)]  We need the  non-positivity of some sectional curvatures of the ambient space for our results to work.

\item[(3)] The geometry of the problem implies that the evolving manifolds have positive averaged mean curvature. This 
 is necessary in proving  the preservation of the  generating curve as a graph in  \cite[Theorem 5]{CaMi2}. Sometimes, this is also a usual restriction asked to get a friendlier  flow behavior (cf. \cite{Li}).
\end{itemize}

\subsection{Suitable ambient spaces}\lb{AS}

Here we give precise definitions of the ambient spaces where we consider the evolution, and also of the concept of revolution hypersurface in them.

\begin{defi} \lb{RSS_def}
A $(n+1)$-dimensional {\bf rotationally symmetric space} (RSS) with respect a curve $\mathcal A$  is a Riemannian manifold $(\oM,\og)$ such that 
there is an  action of $SO(n)$ on $(\oM,\og)$ by isometries for which the set of fixed points is the curve $\mathcal A$. Then $\mathcal A$ is a geodesic and it is called the {\bf axis of rotation}.

A smoothly embedded hypersurface $X:M\flecha \oM$ is said to be a {\bf hypersurface of revolution} around $\mathcal A$ if it is invariant under the action of $SO(n)$ on $(\oM,\og)$.
\end{defi}

There are natural ways of constructing a RSS by using warped products and spherically symmetric spaces. Recall that a  warped product $\mathcal M \times_f \mathcal N$ of two Riemannian manifolds $(\mathcal M, g)$ and $(\mathcal N, h)$ is given by $(\mathcal M\times \mathcal N,  g + f^2 h)$, being $f:\mathcal M\flecha \re$ a positive smooth map. A {\bf spherically symmetric space}  $(\mc S,\sigma)$ admits a metric of the form $\sigma= dr^2 + h(r)^2 g_{{\mathbb S}^{n -1}}$ with $h(0) = 0$ and $h'(0)=1$, where $r$ is the distance to a fixed point $\mathcal O$ in $\mc S$ and $g_{{\mathbb S}^{n -1}}$ is the metric of the round unit sphere. Here we shall consider the more standard complete cases:

$\circ$ \cite[section 3.2]{Gn} $\mathcal O$ is a pole, then  $h$ never vanishes, $\mc S$ is diffeomorphic to $\re^n$ and can be parametrized on $[0,\infty[ \times {\mathbb S}^{n -1}$;  

$\circ$ \cite[page XV.13]{BG} the {\it first positive zero $\z$ of} $h$ exists ($\z<\infty$), then $h(\z)=0$, $h'(\z)=-1$,   $\mc S$ is a differentiable sphere and can be parametrized on $[0,\z[ \times {\mathbb S}^{n -1}$. 

 In short, ${\mc S}$ can be regarded as the warped product $I \times_h {\mathbb S}^{n - 1}$,  with  $I=[0,\z]$ when $\z<\infty$ and $I=[0,\infty[$ otherwise.

In practice, we consider two kinds of warped products to build up a RSS:
$(\tM, \tg):= {\mc S} \times_f J$, $\text{ with } f: {\mc S} \flecha \re$ depending only on $r$ or
$(\oM, \og):= J\times_f \mc S, \text{ with } f:J\flecha \re$ and  $J$ a real interval. The above expression for the metric $\sigma$  yields 
\begin{align}
&(\tM,\  \tg):= (I \times {\mathbb S}^{n -1}\times J,  \  dr^2  +  h(r)^2 g_{{\mathbb S}^{n -1}} + f(r)^2 dz^2)  \lb{meMZ} \\
&\hspace*{-2cm} \text{ and }\nn \\
&(\oM,\  \og):= (J\times I\times {\mathbb S}^{n -1}, \  dz^2 + f(z)^2 dr^2 + f(z)^2 h(r)^2 g_{{\mathbb S}^{n -1}}). \label{omet}
\end {align}
In both cases the action of $SO(n)$ is given by 
$$R(z,r,u) = (z,r,Ru), \text{ for every }R\in SO(n).$$
Obviously $\mc A_+:= J\times\{0\}\times {\mathbb S}^{n -1}$ (with ${\mathbb S}^{n -1}$ collapsed to a point, because $h(0)=0$) is part of the axis of rotation $\mc A$, which coincides with $\mc A_+$ when $\z$ does not exist. If $\z<\infty$, one has $\mc A = \mc A_+ \cup \mc A_-$, with $\mc A_-:= J\times\{\z\}\times {\mathbb S}^{n -1}$ (with ${\mathbb S}^{n -1}$ collapsed to a point, because $h(\z)=0$). 
 
In \cite{CaMi2} we used \eqref{meMZ} as the ambient space. Here shall see (cf. section \ref{RSS}) that the hypesurfaces $z=constant$ in $(\oM,\og)$ are orthogonal to the axis $\mathcal A_+$ and at constant distance from each other. Then \eqref{omet} is specially suited to consider equidistant hypersurfaces as the boundary of the domain containing the surface to evolve. As we shall show in Remark \ref{models}, space forms are special cases of $(\oM, \og)$, and specific choices of the functions $f$ and $h$ give also a new situation in the Euclidean space.

We are thus led to consider the following natural setting:
\begin{quote}
{\bf Setting $\frak{Eq}$}.  $(\oM,\og)$ is a RSS with axis of rotation $\mathcal A$ and metric $\og$ as in \eqref{omet} satisfying either $\int_0^\infty h(r)^{n-1} dr =\infty$ or $\z<\infty$.  $M\subset \oM$ is a smoothly embedded hypersurface of revolution around $\mathcal A$  generated by the graph of a function $r(z)$ over $\mc A_+$ and contained  
in the domain
$G = \{(z,r,u)\in \oM  : a \le z\le b\}$,  with boundary $\partial M$, which intersects $\partial G$ orthogonally  and encloses a $(n+1)$-volume V inside $G$. 
\end{quote}

Then flow $M$  by \eqref{vpmf} with the boundary condition that 
\be \lb{bocon}
\text{$M_t$ intersects $G$ orthogonally at the boundary for every $t$}.
\ee

\subsection{Statement of the main results}\lb{MR}

Throughout this paper we shall prove:

\begin{teor}\label{MT} Let $M_t$ be the solution of  \eqref{vpmf}  with initial condition in the setting $\frak{Eq}$ and boundary condition \eqref{bocon},  defined on a maximal interval $[0,T[$. Then
\begin{itemize}
\item[a)] The generating curve of the solution $M_t$ of \eqref{vpmf} remains a graph over  $\mc A_+$ for every $t\in[0,T[$. 
\item[b)] If $T<\infty$, the singularities at $t=T$  are located on the axis of rotation $\mc A$.
\item[c)] There is a constant $C$ depending on $\og$, $V$, $a$ and $b$ such that if  $|M| \leq C$, then $T=\infty$ and there is a sequence of times $t_n\to\infty$ such that $M_{t_n}$ converges to a revolution hypersurface of constant mean curvature in $\oM$.
\end{itemize}
\end{teor}

This result not only completes the non-Euclidean version of \cite{Ath1,Ath2}, started in \cite{CaMi2}, by considering equidistant instead of totally geodesic hypersurfaces as the boundary of the domain containing the evolution. In fact, it also solves the problem for a new  situation in the Euclidean space: the case where the boundary hypersurfaces are spheres instead of hyperplanes (see Remark \ref{models} for details).

More surprisingly, this change of the boundary hypersurfaces makes the corresponding results valid for a new and interesting framework: ambient spaces with positive curvature and evolving hypersurfaces with non-necessarily positive $\oH$. To our knowledge, besides breaking the restrictions (2) and (3)  of our  statements in \cite{CaMi2}, this is the first time that results for the evolution \eqref{vpmf} are obtained in  a family of ambient spaces of positive cuvature including those of constant curvature, and allowing the possibility $\oH<0$.

Such a new scenario is even more rewarding if we realize that we are in a much harder situation than those in \cite{Ath1, Ath2, CaMi2}. Indeed, the geometry of the new setting does not allow to use any symmetry as we pointed out in (1), so each step in the proof has the further complication of analyzing what happens at the boundary.

The paper is organized as follows. In Section \ref{RSS} we study the geometry of the ambient space with the metric \eqref{omet} and give some special interesting examples of the setting $\frak Eq$. Section \ref{STE} gathers computations of basic quantities for evolving revolution surfaces, standard results about short time existence and basic evolution formulas for our flow. In Section \ref{broH} we obtain upper bounds for the distance to $\mc A_+$ and for the absolute value of the averaged mean curvature, results that we shall apply in Section \ref{graph} to prove the preservation of the property of being a graph for the generating curve of the evolving hypersurface. Section \ref{ie} is devoted to obtain interior  estimates of the heat operator acting on a special function, which is applied in Section \ref{bL} to get more interior estimates, boundary estimates and uniform bounds for the norm of the Weingarten map. In Section \ref{preLTE} we obtain the estimates for the higher order derivatives, concluding with the proof of part b) of Theorem \ref{MT}. Finally, in Section \ref{volp} we prove part c) of the theorem. Appendix A is devoted to the proof of  a computational lemma, and in Appendix B we give two examples of  hypersurfaces in the setting $\frak Eq$ with negative averaged mean curvature.

\section{More about the geometry of the RSS $(\oM, \og)$}\label{RSS}

For subsequent arguments, it will be very useful to have explicit expressions for the Levi-Civita connection $\oN$ of $(\oM,\og)$. Given a local orthonormal frame $\{e_i\}_{i=2}^n$ for the unit sphere $\mathbb S^{n-1}$ with the standard metric and the vector fields $\partial_r, \partial_z$ associated to the coordinates $r$ and $z$ of $\oM$, it follows from the expression of $\og$ that $\{ \partial_z, E_r, E_2, ..., E_n\}$ (with $E_r = \fracc{\partial_r}{f}$, $E_i= \fracc{e_i}{f h}$) is a local orthonormal fame of $(\oM,\og)$. Since $\oM$ is the  warped product $J \times_f \(I \times_h \mathbb S^{n-1}\)$, using the  formulae for the covariant derivatives of a warped product (cf. \cite {ON}), we obtain

\begin{lema} \lb{Conn} For the Levi-Civita connection $\oN$ of $(\oM,\og)$, the following formulae hold
\begin{align}
\oN_{\partial_z} \partial_z &=  0,    \quad  \oN_{\partial_z} E_i   =  0 , \quad  \oN_{\partial_z} E_r  =  0 , \quad \oN_{E_i} \partial_z =   \fracc{f'}{f} E_i ,  \lb{umbz} \\
 \oN_{\partial_z} \partial_r   &= \oN_{\partial_r} \partial_z = \frac{f'}{f} \partial_r , \quad   \quad \oN_{\partial_r} \partial_r  = - f' f \partial_z, \quad   \oN_{E_r} E_r  = - \frac{f'}{f} \partial_z,  \lb{dzr}\\
\oN_{\partial_r}E_i &=  0 , \quad \quad  \quad \oN_{E_i}\partial_r =  \frac{h'}{h} E_i \lb{dzri} \\
  \oN_{E_i} E_j &=   \frac1{(f h)^2}  \nabla^{\mathbb S}_{e_i}e_j -  \(\frac{f'}{f} \partial_z  +  \frac{h'}{f^2h} \partial_r\) \delta_{ij}. \lb{Dij} \end{align}
for $2\le i,j \le n$, where $\nabla^{\mathbb S}$ denotes the Levi-Civita connection of ${\mathbb S}^{n-1}$. 
\end{lema}

\noindent\begin{nota}\lb{Kill} 
It follows from \eqref{dzr} that the \lq\lq plane $zr$" is a totally geodesic surface and that  $\partial_r$ restricted to that surface is a Killing vector field.
\end{nota}

\begin{nota}\lb{Kill2} From \eqref{omet}, \eqref{umbz} and \eqref{dzr} we deduce that the curves $z\mapsto(z, r_0,u_0)$ are geodesics and  the hypersurfaces $z= c$ ($c$ constant) are at constant distance from each other, and umbilical with normal curvature $\fracc{f'}{f}(c)$. 
Hence only the values $c$ of $z$ for which $f'(c)=0$ make the hypersurface $z=c$  totally geodesic in $\oM$. If such a $c$ exists, the boundary hypersurfaces of $G$ in setting $\frak Eq$ are equidistant from a totally geodesic one, which corresponds to a special framework in the Hyperbolic Space (see case  {\rm C3} in Remark \ref{models}).
 
The hypersurfaces $z=$ constant have the same constant normal curvature $k$ if and only if $\fracc{f'}{f}(z) = k$, which gives   $f(z) = d\  e^{k z}$ for some constant $d$. These hypersurfaces correspond to  horospheres when the ambient space is the Hyperbolic Space (case {\rm C4} in Remark \ref{models}).
\end{nota}

Using now the formulae for the curvature of a warped product and the standard expression for the curvature tensor of $\mathbb S^{n-1}$ , we obtain
\begin{lema} \lb{Curv}
The components of the curvature tensor $\oR$ of $(\oM,\og)$ in the basis \linebreak $\{\partial_z, E_r, E_2, ..., E_n\}$ are 
\begin{align*}
\oR_{z \a  \beta \gamma}  & = 0 \qquad \oR_{z\a z \beta}  = -\fracc{f''}{f} \delta_{\a \beta}   \qquad  \\
\oR_{r i j k} &= 0 \qquad \oR_{r i r j}  = - \frac1{f^2}\( \frac{h''}{h} + f'^2\) \delta_{i j} \qquad \\
\oR_{i j k \ell}  
&=  \fracc{1-(h'^2 + h^2 f'^2)}{f^2 h^2}  \(\delta_{ki} \delta_{\ell j} - \delta_{\ell i} \delta_{kj}\) \qquad
\end{align*}
for $\a, \beta, \gamma  \in\{r, 2, \ldots, n\}$ and $i, j, k, \ell \in \{2, \ldots, n\}$.
\end{lema}

\begin{nota} \lb{models}
From lemmas \ref{Conn} and \ref{Curv} we have the following different special cases for the setting $\frak Eq$ in space forms: 
\begin{itemize}
\item[\bf(C1)] $J=\re$, $f(z)= 1$ and $h(r)=r$. Then $I =[0,\infty[$, $(\oM,\og)$ is the Euclidean space $\re^{n+1}$ and $G$ is the slice between two hyperplanes. $\mc A = \mc A_+ $  is the axis $x^{n + 1}$ in $\re^{n+1}$.
\item[\bf(C2)] $J=[0,\infty[$, $f(z)= z$ and $h(r)=\sin r$. Then $I =[0,\pi]$, $(\oM,\og)$ is  again the Euclidean space  and $G$ is the spherical crown between two spheres of radii $a$ and $b$. $ \mc A_+$ is the upper half-axis $x^{n+1}$ in $\re^{n+1}$ and $ \mc A_-$ is the lower half-axis.
\item[\bf (C3)] $J=\re$, $f(z)=\cosh (\sqrt {|\lambda|}\  z)$ and $h(r)=  {|\lambda|}^{-\frac1{2}} \sinh (\sqrt {|\lambda|}\  r)$ for $\lambda <0$. Then $I =[0,\infty[$, $(\oM,\og)$ is ${\mathbb H}^{n +1}_\lambda$, which here means  the Hyperbolic space of sectional curvature $\lambda$,  and $G$ is the slice between two equidistant hypersurfaces. $\mc A = \mc A_+ $. 
\item[\bf (C4)] $J=[0,\infty[$,  $f(z)=   {|\lambda|}^{-\frac1{2}} \sinh (\sqrt {|\lambda|}\  z)$ and $h(r)=  \sin\ r$ for $\lambda <0$. Then $I =[0,\pi]$, $(\oM,\og)$ is ${\mathbb H}^{n + 1}_\lambda$  and $G$  is the spherical crown between two geodesic spheres of $\oM$ of radii $a$ and $b$.
\item[\bf (C5)] $J=\re$, $f(z)=e^{\ral z}$ and $h(r)= r$, $\lambda <0$. Then $I =[0,\infty[ $, $(\oM,\og)$ recovers again  ${\mathbb H}^{n + 1}_\lambda$ and $G$ is the slice between two \lq\lq parallel'' horospheres. $\mc A = \mc A_+ $.
\item[\bf (C6)]$ J=\left[-\fracc{\pi}{2 \sqrt{\lambda}}, \fracc{\pi}{2 \sqrt {\lambda}}\right]$,  $f(z)=\cos (\sqrt {\lambda}\  z)$ and $h(r)=  {\lambda}^{-\frac1{2}} \sin (\sqrt {\lambda}\  r)$ for $\lambda >0$. Then $I =\left[0, \pi/\sqrt{\lambda}\right]$, $(\oM,\og)$ is the round ${\mathbb S}^{n +1}(1/\sqrt{\lambda})$ and $G$ is the slice between two parallels. $ \mc A=\mc A_+ \cup \mc A_-$ is a meridian, with $\mc A_+$ and $\mc A_-$ half-meridians.
\end{itemize}

Let us remark that even in the cases (C2), (C4) and (C6) where $\mc A\ne \mc A_+$, one has that  $\mc A=\mc A_+ \cup \mc A_-$ is a connected real line (a circle in case (C6)); accordingly,  even   from an intuitive viewpoint, $\mc A$ has the right to be called the axis of rotation.
\end{nota}
\begin{nota}\lb{models2}
If in examples (C3) and (C6) we use different constants in the definition of $f$ and $h$ (for instance, in (C3) we pick $f =\cosh(\sqrt{\lambda} \, z)$, $h =  {|\mu|}^{-\frac1{2}} \sinh (\sqrt {|\mu|}\  r)$ with $\mu \neq \lambda$), we still produce constant sectional curvature, but we get spaces with singularities (or not complete regular spaces). These model ambient spaces appear in the literature as extremals of some functionals defined on the space of Riemannian metrics (cf. \cite{Ks, GP, E, LlM}). Since our theorem refers to slices $G$ which do not contain the singular points, it is also true in these non-regular ambient spaces. 
\end{nota}

\section{Evolving revolution hypersurfaces within a RSS}\label{STE}

Let us begin with a remark on the notation: when  we introduce for the first time a quantity depending on the evolving hypersurface $M_t$, we write either a subindex $_t$ or $ (\, .\, , t)$ to denote its dependence on $t$, but just later we shall omit this unless it is not clear from the context what we mean.

\subsection{Basic quantities on revolution hypersurfaces}

Our flow \eqref{vpmf} is invariant under isometries of $(\oM,\og)$, then it is invariant under the action of $SO(n)$. As a consequence, if the starting hypersurface $M$ is of revolution, also is the evolving $M_t$. Hence  the unit normal vector $N_t$ to $M_t$ will be contained in the plane generated by $E_r$, $\partial_z$ and can be written as
\be
N= \<N,E_r\> E_r + \<N,\partial_z\> \partial_z, \label{N}
\ee
in turn, the unit vector $\frak t_t$ tangent to the generating curve will be
\be
\frak t = -\<N, \partial_z\> E_r + \<N, E_r\> \partial_z. \label{t}
\ee

We shall use the coordinates $(z_t,r_t,u_t)$ for $M_t$. 
Without loss of generality, we can parametrize the generating curve $c_t$ of $M_t$ as $c: [a,b] \flecha \oM$, $c(s) = (z(s), r(s), u)$, with $\dot{c}(s)\ne 0$ for every $s$. With this parametrization, the vectors $\frak t$ and $N$ admit the expressions
\be\label{tNf}
\frak t = \frac{1}{|\dot{c}|}(f \dot{r} E_r +  \dot{z} \partial_z), \quad N= \frac{1}{|\dot{c}|} (\dot{z}  E_r - f \dot{r} \partial_z),
\ee
where $|\dot{c}| := \sqrt{\dot{z}^2+ (\dot{r}\, f)^2}$ and  $\dot{z}$, $\dot{r}$ denote the derivatives of $z$ and  $r$ with respect to $s$. 

Consider the local orthonormal frame  $\frak t, E_2, ... , E_n$ on $M_t$. Then the mean curvature of $M_t$ is given by
\be\label{Hf1} H = k_1 + (n-1) k_2, \ee
where $k_1$  is the normal curvature of $M_t$ in the direction of $\frak t$:
\begin{align} 
k_1 = - \<\oN_{{\frak t}} {\frak t}, {N}\>
&= - \frac{1}{|\dot{c}| } \left( \frac{\ddot{r} f \dot{z} - \ddot{z} f \dot{r}  + \dot{r} f' \dot{z}^2}{|\dot{c}| ^2} + f' \dot{r} \right), \label{k1f1}
\end{align}
and $k_2$ is the normal curvature of $M_t$ in the direction of $E_i,\ i=2, ..., n$:
\begin{align}\label{k2f2} k_2 &= \<\oN_{E_2} N, E_2\> = \<N, E_r\> \<\oN_{E_2} E_r, E_2\> + \<N, \partial_z\> \<\oN_{E_2} \partial_z, E_2\> \nn
\\ & = \frac{h'}{h f} \<N, E_r\> + \frac{f'}{f} \<N, \partial_z\> =  \frac{1}{|\dot{c}|} \( \frac{h' \dot{z} }{h f} - f' \dot{r} \).
\end{align}

\subsection{Short time existence and some evolution formulae}
Recall the well known fact (cf. \cite{Eck}) that  $X_t$ is a solution of \eqref{vpmf}  if and only if it is, up to tangential diffeomorphisms, a solution of 
\be\label{vpmft}
\<\parcial {X}{t}, N\> = \oH - H.
\ee

If we consider the flow of the graph of $(z,u)\mapsto (z,r(z),u)$ under \eqref{vpmft}, the variable $z$ does not change with time, and formulae of the previous subsection (now taking $s=z$) remain true for any time. Using them,  equation \eqref{vpmft} with this initial condition becomes
\begin{align}
\parcial{r}{t} &=   \frac{ \ddot{r} }{|\dot{c}| ^2} + \frac{f'}{f}\( \frac{1}{|\dot{c}| ^2} +n\) \dot{r}- (n-1)  \frac{h'}{h f^2}+ \oH \frac{|\dot{c}| }{f}.\label{tafl}
\end{align}
Here replacing $\oH$ in \eqref{tafl} by any $C^{1, \alpha/2}$ real valued function $\psi$ such that $\psi(0) = \oH(0)$, we obtain a parabolic equation which, at least for small $t$, has a unique solution satisfying $\dot{r}(a) =\dot{r}(b)=0$. Now, using a routine fixed point argument (cf. \cite{Mc3}), we can establish short time existence also for \eqref{tafl} with the same boundary conditions.

\medskip

The following lemma collects some  evolution formulae for \eqref{vpmf} in $(\oM,\og)$.
\begin{lema} \lb{Evol_gen}
If $M_t$ is a solution of \eqref{vpmf}, the following evolution equations hold:
\begin{align*}
\text{\rm (a)} & \ \ds \frac{\oN N}{\partial t} =  \nabla H  \\
\text{\rm (b)} &\ \ds \deri{}{t} |M_t| =  - \int_M (\oH - H)^2 \ d\mu_t,  \\
\text{\rm (c)} &\ \parcial{|L|^2}{t} = \Delta |L|^2 - 2|\nabla L|^2 + 2|L|^4- 2 \oH  \vtr L^3 +  2 |L|^2 \(\mc{T} + (n-1) \mc{J}\)  \nn \\
& \quad  - 2 \oH \(k_1 \mc{T} + (n-1) k_2 \mc{J} \) - 4 (k_1- k_2)^2 (n-1) \mc{Y} - 2 \< \a  , \tilde{\delta} \oR_N \>, 
\end{align*}
where $\nabla$  denotes both the intrinsic covariant derivative and the gradient on $M_t$, $\Delta$ denotes its intrinsic Laplacian and $\a$ its second fundamental form. Moreover  $\tilde{\delta} \oR_N(X,Y):= \sum_i\(\nablaa_X \oR_{N E_i Y E_i} + \nablaa_ {E_i}\oR_{N Y X E_i} \)$, $ \mc{T}$, $\mc{J} $ and $\mc{Y}$ are the sectional curvatures of the planes generated by $\{\vt, N\}$, $\{E_i, N\}$ and $\{E_i, \vt\}$, respectively. 
\end{lema}
\begin{demo}
(a) and (b) are well known and valid for any ambient space. The proof of (c) follows arguing exactly the same that in \cite{CaMi2}, substituting the orthonormal local frame $N,\frak t=E_1, E_2, ..., E_n$ into the more general  and standard evolution equation of $|L|^2$ (see, for instance, (6.1) in \cite{CaMi2}).
\end{demo}

The lemma below contains two equations which are very specific to our setting. The proof is straightforward but quite long and technical; the interested reader can find the details in the Appendix  A of the present paper.
\begin{lema} \lb{Evolrz}
Set $u := \<N, \partial_r\>$, then for any functions $\phi, \psi:\re\flecha \re$ one has the following evolution formulae under \eqref{vpmf}:

\medskip 
{\rm (a)} $\ds \(\parcial{}{t}-\Delta\) \phi(r) = \phi'\left( \oH \frac{u}{f^2} - 2 \frac{f'}{f^3} u \<N, \partial_z\> - (n - 1) \frac{h'}{f^2 h} \right) +  \frac{\phi''}{f^2} \(\frac{u^2}{f^2} - 1\).$

\medskip 
{\rm (b)} $\ds \(\parcial{}{t}-\Delta\)\psi(z) =  \psi' \(\oH  \<N,\partial_z\> + \(\frac{u^2}{ f^2} - n\) \frac{f'}{f} \) - \psi'' \frac{u^2}{f^2}.$
\end{lema}

\section{Upper bounds for $r$ and  $|\oH|$}\label{broH}
 
In this section, we shall prove  that if $M$ is a hypersurface satisfying the conditions in the setting $\frak{Eq}$, then the coordinate $r$ of $M$ and its averaged mean curvature $\oH$ (in absolute value) are uniformly bounded. In fact, we shall obtain these bounds under more general conditions than the setting $\frak Eq$. Our ultimate goal is to apply these results to bound $r$ and $\oH$ for  a maximal solution $M_t$  of \eqref{vpmf}. 

\begin{notacion}\lb{infinito} 
From now on, given any function $F(z,r,u)$, we shall use the notation $\| F\|_\infty = \sup_{[a,b]\times[\rho,\de]\times \Sn} |F(z,r,u)|$, where $\rho$ and $\de$ are constants explained in each situation and many times $F$ depends only on one or two of the variables $z$, $r$, $u$.
\end{notacion}

Let us define the function $\delta(R)=\ds\int_0^R  h(r)^{n-1} dr,$ and let $r_2$  be the constant 
\be
 r_2 =\delta^{-1}\left(\fracc{|M|}{\omega_{n-1}} \ \|f^{-n}\|_\infty + \fracc{V}{\omega_{n-1} \int_a^b f(z)^n dz}\right),  \label{defr2}
\ee
where  $\omega_{n - 1}$ denotes the volume of $\mathbb S^{n-1}$ with its standard metric. When $\z=\infty$ (which hereafter means that $\z<\infty$ is false), the hypothesis $\int_0^\infty h(r)^{n-1} dr = \infty$ in  setting $\frak Eq$ ensures that $r_2$ always exists. By contrary, when $\z<\infty$, $r_2$ may not be well defined; if this happens, we use the convention  $\min\{\z,r_2\} = \z$.

Observe that the following result does not require the generating curve of $M$ to be the graph  of a function nor contained in $G$.

\begin{prop}\label{borM}  Let $(\oM,\og)$ be as defined in \eqref{omet}. If $M$ is an embedded hypersurface of revolution in $\oM$, with boundary in the hypersurfaces   $z=a$, $z=b$ and orthogonal to them along the boundary, 
 then
$r< \min\{ \z, r_2\}$.

\end{prop}

\begin{demo} 
If there is some point in $M$ with $r = \z$, as $M$ is of revolution, this point has to be singular, in contradiction with the fact that $M$ is a regular submanifold. Then, we shall concentrate on proving that $r<r_2$, with $r_2 < \z$. 

Now we define 
\be \lb{defR1}
r_1 = \delta^{-1}\(\fracc{V}{\omega_{n-1} \int_a^b f(z)^n dz}\).
\ee
 It follows that $r_2> r_1 >0$ because $\delta$ is an increasing function, 
 
 Let us denote by $r_m$ and $r_M$ the minimum and maximum value of $r$ on $M$ respectively, and let $r_z= \inf\{r(s); z(s)=z\}$. By  $d \mu_{\og}$ we mean the volume element of $\oM$ and by $\Omega$ the domain enclosed by $M$ and the disks  in $\partial G$ limited by $\partial M$. Using  the definition of $r_1$ and the expression \eqref{omet}, we obtain
\begin{align}
 \omega_{n - 1} & \( \int_a^b f(z)^n dz\) \(\int_0^{r_1}  h(r)^{n-1} dr\) = V= \int_{\Omega} d \mu_{\og}  \nn\\
 & \ge \omega_{n - 1} \int_a^b \int_0^{r_z} f(z)^n h(r)^{n-1} dr \ dz  \nn\\
 & \ge  \omega_{n - 1}\( \int_a^b f(z)^n dz\) \(\int_0^{r_m}  h(r)^{n-1} dr\) ,  \label{fVr1}
\end{align}
Next, recalling \eqref{omet}, \eqref{tNf} and that $d\mu = \iota_N d\mu_\og$, we get the area of $M$ as
\begin{align}
|M| &= \int_M  \imath_{N} \(d \mu_{\bar g}\) = \omega_{n - 1} \int_a^b \sqrt{\dot{z}(s)^2 + f(z(s))^2 \dot{r}(s)^2 } \  f(z(s))^{n-1} h(r(s))^{n-1} ds\label{muf}\\
&> \omega_{n - 1} \int_a^b   |\dot{r}(s)| \  f(z(s))^{n} h(r(s))^{n-1} ds \ge  \omega_{n - 1} \min_{[a,b]} f(z)^n  \int_{r_m}^{r_M}    h(r)^{n-1} dr \label{vle>+}.
\end{align}
From the inequality \eqref{fVr1} we have that $r_1 \ge r_m$. If $r_1 \ge r_M$, we have the desired bound. If not, it follows from the inequality \eqref{vle>+} that 
\begin{align}
|M| & >   \frac{\omega_{n - 1}}{\|f^{-n}\|_\infty}  \left[\int_{r_m}^{r_1} + \int_{r_1}^{r_M} \right]    h(r)^{n-1} dr  > \frac{\omega_{n - 1}}{\|f^{-n}\|_\infty}    \int_{r_1}^{r_M}     h(r)^{n-1} dr \nn \\
& = \frac{\omega_{n - 1}}{\|f^{-n}\|_\infty}  \  (\delta(r_M) - \delta(r_1)).\nn
\end{align}
Hence
\be\lb{delrM}
\delta(r_M) < \fracc{|M| \|f^{-n}\|_\infty}{ \omega_{n - 1} } + \delta(r_1),
\ee
from which the proposition follows. 
\end{demo}

\begin{coro}\label{bor}
If $(\oM,\og)$ and $M$ are in the setting $\frak{Eq}$  and $[0,T[$ is the maximal time interval  where the flow  \eqref{vpmf} satisfying the boundary condition \eqref{bocon} is defined; then
$r_t <\min\{\z, r_2\}$ for every $t\in[0,T[$, with $r_2$ defined by \eqref{defr2} for the initial condition $M_0$.
\end{coro}

\begin{demo} Applying Proposition \ref{borM} to $M_t$ for each fixed $t$, we reach \eqref{delrM} with $|M_t|$ instead of $|M|$.  Then the conclusion follows using the area decreasing property of the flow  (which is a consequence of (b) in Lemma \ref{Evol_gen}) and that  the function $\delta$ is increasing. 
\end{demo}

Next, the goal is to bound the modulus of the averaged mean curvature $\oH_t$.

\begin{prop}\label{p:uboHM}  Let $(\oM,\og)$ and  $M$ as in Proposition \ref{borM}. If the number of points in the  generating curve of $M$ with tangent in the direction of $\partial_r$ is finite and $\z > \de \ge r\ge \rho>0$,
 then there is a constant $h_2(V, \overline g,n,a,b,\rho, \de)> 0$   such that $| \oH | \le h_2$.
\end{prop}

\begin{nota}  Observe that when $r_2 < \z$ the hypothesis $\de \ge r$ suppose no restriction, since by Proposition \ref{borM} we can take $\de = r_2$. In this case $h_2$ depends on $|M|$ through $r_2$.
\end{nota}

\begin{demo} 
From \eqref{Hf1}, \eqref{k1f1} and \eqref{k2f2} we can write
\begin{align}
\oH = \frac1{|M|} \int_{M}  \frac{-1}{|\dot{c}| }  \frac{d}{ds}\arctan\left({f} \frac{\dot{r}}{\dot{z}}\right) d\mu +  \frac1{|M|} \int_{M} \frac{1}{|\dot{c}| }\left((n-1)\frac{h'}{h f} \dot{z} - n f' \dot{r}\right) d\mu =: I_1 + I_2. \label{oHt}
\end{align}

Now we integrate by parts,  and having into account that  the condition in the boundary gives $\dot{r}(b)=\dot{r}(a)=0$ and that at the points $s_i$, $i=1, \dots, k$ where the tangent vector  to the generating curve is vertical (that is $\fracc{\dot{r}}{\dot{z}}=\pm \infty$) one  still has that $\arctan\( f(z(s_i)) \fracc{\dot{r}(s_i)}{\dot{z}(s_i)} \)$ is finite,  we get
\begin{align}
I_1 &= \frac{\omega_{n - 1}}{|M|}  \int_a^b   - (f h)^{n-1} \frac{d}{ds}\arctan\left(f\, \frac{\dot{r}}{\dot{z}}\right)ds 
\quad   \nn \\
& =\frac{(n-1) \omega_{n - 1}}{|M|} \int_a^b \arctan\left(  f \frac{\dot{r}}{\dot{z}} \right) (fh)^{n-2}\  (f h' \dot{r} + h f' \dot{z})\  ds. \label{I1oH} 
\end{align}
Using $\arctan\left(f \fracc{\dot{r}}{\dot{z}} \right) f \fracc{\dot{r}}{\dot{z}}  \le \fracc{\pi}{2 |\dot{z}|} |f \dot{r}|$, $|f \dot{r}|\le |\dot{c}|$ and $| \dot{z}|\le |\dot{c}|$, we get 
\begin{align}
|I_1|  &  < \frac{(n-1)\, \omega_{n - 1}}{|M|} \frac{\pi}{2} \( \int_a^b  |\dot{c}| (f h)^{n -2}  |h'| \, ds + \int_a^b f^{n-2} h^{n-1} |f'| |\dot{c}|\,ds\) \nn \\
&\le \frac{(n-1)}{|M|} \frac{\pi}{2}  \( \int_M \frac{|h'|}{f h} d\mu + \int_M  \frac{|f'|}{f} d\mu\)  \leq (n -1)\frac{\pi}{2} \(\left\| \frac{h'}{hf}\right\|_{_\infty}+  \left\| \frac{f'}{f} \right\|_{_\infty}\). 
\label{I1ub}
\end{align}

Next, we bound $|I_2|$ as follows:
\be
|I_2 |\le (n-1) \left\|\fracc{h'}{fh}\right\|_{_\infty} +  \frac{n}{|M|} \int_M \fracc{ |f'\dot{r}|}{|\dot{c}|} \, d\mu  <  (n-1) \left\|\fracc{h'}{fh}\right\|_{_\infty} + n  \left\|\fracc{f'}{f}\right\|_{_\infty} \label{I2pos}
\ee
In conclusion, the existence of the finite upper bound $h_2$ follows from \eqref{I1ub} and \eqref{I2pos}. 
\end{demo}

\begin{coro}\label{p:uboH}
Let $M_t$ be the solution of  \eqref{vpmf}  with initial condition $M$ in the setting $\frak{Eq}$ and satisfying the boundary condition \eqref{bocon}.
For every $t$ such that $0<\rho \le r_t \le \de < \z$ and the generating curve of $M_t$ is a graph, 
there is a constant $h_2(V, \overline g,n,a,b,\rho,\de)> 0$   such that $| \oH | \le h_2$.
\end{coro}

\begin{demo} It follows because if the generating curve of $M$ is a graph, it satisfies the conditions  in Proposition \ref{p:uboHM}.
\end{demo}

\section{The generating curve remains a graph}\label{graph}

This section is devoted to prove that, for $M$ in the setting ${\frak Eq}$, the evolving hypersurface remains  a revolution hypersurface generated by a smooth graph. As we pointed out above, $M_t$ is always a revolution hypersurface. Then the aim is to show that the generating curve remains  a graph over the axis of rotation for all time.

Recall that the generating curve is a graph if and only if $\< N, E_r\> >0$, which is equivalent to say $u:= \<N,\partial_r\> >0$, and also equivalent to $\fracc{1}{f(z)}\le v = \fracc{1}{u}< \infty$. Therefore, our goal is to obtain an upper bound for $v$. To achieve this,  we first need  the evolution equation for $v$.

\begin{lema} \label{dvtn} Under \eqref{vpmf}, $v= u^{-1}$ evolves as
\begin{align*} 
\parcial{}{t} v &= \Delta v -  \(|L|^2 + \oRic(N,N) + \frac{(n - 1)}{f^2} \left[\frac{h'}{h}\right]' \) v - \frac2v |\nabla v|^2.
\end{align*}
\end{lema}

\begin{demo} 
First we compute $\Delta u$. To do so, for a fixed time $t$ and a point $p\in M$, we shall use a local  frame $F_1, F_2, ... , F_n$ of $M$ orthonormal at $p$ and satisfying $\oN_{F_i}F_j(p)=0$. Next, we extend it to a local frame $\gor F_i$ on a neighborhood of $p$ in $\oM$ using the flow of $\partial_r$ so that $[\partial_r, \gor F_i]=0$. It follows from Bartnik's formula (cf. \cite{Bar} page 158) that, at point $p$, 
\begin{align}
\Delta u = & \<\partial_r,  \nabla H\> -(|L|^2 + \oRic(N,N)) u \label{Deltau} \\
& + \sum_{i = 1}^n  \left[\(\oN_{ \gor F_i }\mathcal L_{\partial_r} g\)(N, \gor F_i) - \fracc{1}2 \(\oN_N \mathcal L_{\partial_r} g\)(\gor F_i,\gor F_i)\right]   +\<\mathcal L_{\partial_r} g, \, \alpha\> - \fracc{H}{2} (\mathcal L_{\partial_r} g)(N,N),\nn
\end{align}
where the last term vanishes by Remark \ref{Kill}.

Since  \eqref{Deltau} is evaluated at a point $p$ where $\{\gor F_i\}$ is orthonormal and the relevant expressions  are tensorial, we can use henceforth the frame $\{\frak t, E_2, ..., E_n\}$ (which satisfies $L E_i = k_2 E_i$) instead of $\{\gor F_i\}$. Doing so, using  Remark \ref{Kill} (which also implies $\oN_N \vt \in{\rm span}\{\partial_z, \partial_r\}$) and  \eqref{dzri}, we get
\begin{align}
\sum_{i = 1}^n &\(\oN_N \mathcal L_{\partial_r} g\)(\gor F_i,\gor F_i) =  \sum_{i=2}^n \(\oN_N \mathcal L_{\partial_r} g\)(E_i,E_i) \nn \\
&= \sum_{i=2}^n \Big(2 N \big<\! \oN_{E_i}{\partial_r},E_i\big> - 2 \(\mathcal L_{\partial_r} g\) (\oN_N E_i , E_i) \Big) =  \sum_{i=2}^n 2 N\<\frac{h'}{h} E_i , E_i\> \nn \\
&  = 2 (n - 1) \left[\frac{h'}{h}\right]'\<\oN r, N\> = 2 (n - 1) \left[\frac{h'}{h}\right]' \frac{u}{f^2}. \label{FormA}
\end{align}
Here we have used that 
\be \label{gradr}
\oN r =  \frac1{f} E_r =  \frac1{f^2} \partial_r. 
\ee

Next, exploiting once more Remark \ref{Kill} (which yields $\oN_\vt \vt, \oN_\vt N  \in{\rm span}\{\partial_z, \partial_r\}$), together with \eqref{Dij}, \eqref{dzr} and \eqref{dzri}, we compute
\begin{align}
\sum_{i = 1}^n & \(\oN_{ \gor F_i }\mathcal L_{\partial_r} g\)(N, \gor F_i) + \<\mathcal L_{\partial_r} g, \, \alpha\> \nn \\
& = \sum_{i=2}^n \(\oN_{E_i} \mathcal L_{\partial_r} g\)(E_i, N) + (\mathcal L_{\partial_r} g)(\frak t,\frak t) \a(\frak t, \frak t) + \sum_{i=2}^n (\mathcal L_{\partial_r} g)(E_i,E_i)  \ \a(E_i,E_i) =
\nn \\
&= \sum_{i=2}^n  \Big(E_i\left[\mc L_{\partial_r} g (E_i, N)\right] -  \mathcal L_{\partial_r} g (\oN_{E_i}E_i, N) -   \mathcal L_{\partial_r} g (E_i, \oN_{E_i}N)  + \mathcal L_{\partial_r} g(k_i E_i, E_i)\Big)\nn
\\ & = \sum_{i=2}^n E_i \Big(  \< \oN_{N}{\partial_r} ,E_i\>+ \< \oN_{E_i}{\partial_r} ,N\>\Big) = 0. 
\lb{FormB}\end{align}

Now, plugging \eqref{FormA} and \eqref{FormB} into \eqref{Deltau}, we conclude
\begin{align}
\Delta u &= \<\partial_r,  \nabla H\> -(|L|^2 + \oRic(N,N)) u 
- (n - 1) \left[\frac{h'}{h}\right]' \frac{u}{f^2}. \lb{Lapu}
\end{align}

On the other hand, using part (a) of Lemma \ref{Evol_gen}, \eqref{dzr} and the flow equation \eqref{vpmf},  we get
\begin{align*}
\parcial{}{t} u &= \<\oN_{\partial_t} N, \partial_r\> + \<N, \oN_{\partial_t} \partial_r\>
\nn \\
&= \<\nabla H, \partial_r\> + (\oH-H)  \(\<N, E_r\> \big<N, \oN_{E_r} \partial_r \big> + \<N, \partial_z\> \big<N, \oN_{\partial_z} \partial_r\big>\) 
& = \<\nabla H, \partial_r\> .
\end{align*}
 The substitution of \eqref{Lapu} in the above formula yields
\begin{align} \nn
\parcial{}{t} u &= \Delta u + (|L|^2 + \oRic(N,N)) u  + (n - 1) \left[\frac{h'}{h}\right]' \frac{u}{f^2}, \end{align}
which joint with  the transformation formulae
$$dv = - \frac1{u^2} du, \quad \parcial{v}{t} = - \frac1{u^2} \parcial{u}{t}, \quad \Delta v = -\frac1{u^2} \Delta u + \frac{2}{u^3} |du|^2$$
lead to the equality in the statement.
\end{demo}

Notice that, unlike the corresponding situation in \cite{CaMi2}, we cannot use directly the evolution equation from Lemma \ref{dvtn} in a maximum principle argument to deduce the sought bound for $v$. Instead of $v$, we need to argue with its product by an appropriate function of $r$, as can be seen in the following proof.

\begin{teor}\label{graphbo} 
Let $M_t$ be the solution of  \eqref{vpmf}  defined on a maximal time interval $[0,T[$, with initial condition $M$ in the setting ${\frak Eq}$ and satisfying the boundary condition \eqref{bocon}. Then the generating curve of the solution $M_t$ of \eqref{vpmf} remains a graph over the axis of revolution for every $t\in[0,T[$. 
\end{teor}

\begin{demo} Let us define $\Phi =\phi(r) v$ for some $\phi:\re \flecha \re$. Thanks to part (a) of Lemma \ref{Evolrz} and Lemma \ref{dvtn}, we obtain 
\begin{align*}
\(\parcial{}{t}-\Delta\) \Phi & =v \(\parcial{}{t} - \Delta\) \phi + \phi\(\parcial{}{t}-\Delta\) v - 2 \<\nabla \phi, \nabla v\>   \\
&= \phi' \,  \(\frac{\oH}{f^2} - 2 \frac{f'}{f^3}  \<N, \partial_z\> - (n -1) \frac{h' v}{h f^2} \) +  \frac{\phi''}{f^2}\(\frac1{f^2 v} - v\)  \\
& \quad - \phi \, v\(|L|^2 + \oRic(N,N) +  \frac{(n-1)}{f^2}\left[\frac{h'}{h}\right]'\)  - 2 \frac{\phi}{v} |\nabla v|^2  - 2 \<\nabla \phi, \nabla v\>. 
\end{align*}
Using $- v\<\nabla\phi, \nabla v\> = - \<\nabla \Phi, \nabla v\> + \phi |\nabla v|^2$ and neglecting the term with $|L|^2$, we reach the inequality 
\begin{align*}
& \(\parcial{}{t}-\Delta\) \Phi  
 \leq \frac{\phi'}{f^2} \(\oH - 2 \frac{f'}{f} \<N, \partial_z\>\) + \frac{\phi''}{f^4 v} - \frac{2}{v}\<\nabla \Phi,  \nabla v\> - \frac{\Phi}{f^2} \ \frak{A},
\\ & \text{ with } \ \frak{A} = \oRic(N,N) f^2 + (n - 1)\left[\frac{h'}{h}\right]' + (n - 1)\frac{h' \phi'}{h \phi} + \frac{\phi''}{\phi}.\end{align*}

On the other hand, given any $t_0\in[0,T[$, we have $\min_{M \times [0,t_0]} r(\cdot , t)= \rho(t_0)>0$ and   $\max_{M \times [0,t_0]} r(\cdot , t) = \de(t_0) < \min\{\z, r_2\}$ (thanks to Corollary \ref{bor} and because if $r$ attains the values $0$ or $\z$ at some point in some time, the solution $M_t$ has a  singularity at this point and time).  Now choose $t_1$ as the maximum time in $[0, t_0]$  such that the generating curve of $M_t$ is a graph for every $t\in [0, t_1[$. 

If we take $\phi(r) :=  e^{C r}$, it holds that 
$\phi' = C \, \phi$ and $\phi'' = C^2 \phi$. Accordingly, 
\begin{align*}
\frak{A} & \geq  - \left|\oRic(N,N) f^2 + (n -1) \left[\frac{h'}{h}\right]'\right| + C \left[(n -1) \frac{h'}{h} + C\right] 
\\ &  \geq - \frak R + C\(C - (n -1) \left\|
h'/h\right\|_{_\infty}\),
\end{align*} 
\be \lb{defR}
\text{where } \  \frak R:= \|f^2\|_{_\infty} \| \oRic \|_\infty +  ( n -1) \( \|h''/h\|_{_\infty}  + \|h'^2/h^2\|_{_\infty}\),
\ee

Next,  we can define the constant
$C : =  \frak R + (n -1) \left\|
h'/h\right\|_{_\infty} + 1 < \infty,$
which (since $C \geq 1$) yields 
\begin{align*}
\frak{A} & \geq - \frak R + C( \frak R + 1) \geq C >0.
\end{align*}

Then,  applying Corollary \ref{p:uboH} on $[0, t_1[$, we reach
\begin{align}
\(\parcial{}{t}-\Delta\) &\Phi \le \frak h - \frac{2}{v} \<\nabla \Phi, \nabla v\>   - \widetilde C\   \Phi, 
\end{align}
with 
 $\frak h:= \frak h(V, \og, n, a,b,  \rho(t_0), \de(t_0)) = C e^{C r_2} \|f^{-2}\|_\infty  \(h_2 + 2 \left\| f'/f \right\|_{_\infty} + C \| f^{-1}\|_\infty\)$ and $\widetilde C = C/\|f^2\|_\infty$.
From here,  by application of the maximum principle, we conclude 
\be \lb{Bov}
v \leq  e^{Cr} v = \Phi \le \max\{e^{C \de(t_0)} \max_{M_0} v, {\frak h}/{\widetilde C}\}   \text{ on }[0,t_1[.
\ee
Since the solution $M_t$ is defined, and continuous in $t$, on $[0, T[ \supset [0, t_0[ \supset [0, t_1[$ the bound \eqref{Bov} is true on the whole the interval $[0, t_1]$. Then by continuity of $v$, $v$ will be still bounded on $[0, t_1 + \eps[$, in contradiction with the definition of $t_1$ if $t_1 < t_0$. In conclusion, $t_1=t_0$ and the generating curve of $M_t$ is a graph along all $[0, t_0]$. As $t_0$ is arbitrary, this is true for $[0, T[$.
\end{demo}

\section{Preliminary interior  estimates}\lb{ie}

Here we begin our way towards getting global estimates of $|L|$. Following \cite{CaMi2}, we start by obtaining an interior estimate for the heat operator acting on a certain function of the form $\p(v) |L|^2$. 

\begin{lema}\lb{evolg} Let $M_t$ be a solution of \eqref{vpmf} defined on $[0,T[$ with initial condition $M$ in the setting $\frak Eq$  satisfying the boundary condition \eqref{bocon} and such that there are constants $\de$ and $\rho$ satisfying $\z > \de \ge r_t \geq \rho > 0$ on $[0, T[$. Let ${\frak g} = \varphi(v) |L|^2$, where $\varphi$ is defined by
\be\label{fi}
\varphi(v) := \frac{v^2}{1-k v^2} \qquad \text{with }\qquad k:= \frac1{2\max{v^2}}.
\ee
Then we can find two positive constants $K_1$, $K_2$ so that 
\begin{align*}
\(\parcial{}{t} -  \Delta\) {\frak g} \le - k {\frak g}^2 +
K_1 {\frak g} + K_2 \sqrt{{\frak g}}  - \frac1{\varphi} \< \nabla{\frak g} , \nabla\varphi\> + \frac{v^2 - \p}{2 \p^3} |\nabla \p|^2 \frak g.
\end{align*}
\end{lema}

 Notice that $k$ is a well defined constant depending on $V, \og, n, a, b, \rho$ and $\de$, as follows from \eqref{Bov}.

\begin{demo}
The evolution of ${\frak g}$ is given by those of  $\varphi$ and $|L|^2$ according to the formula:
\begin{align*}
\(\parcial{}{t} -  \Delta\) {\frak g} &= |L|^2 \(\parcial{}{t} - \Delta\) \varphi +  \varphi \(\parcial{}{t} - \Delta\) |L|^2 - 2 \<\nabla \varphi, \nabla |L|^2\> \\ 
& \leq \p'|L|^2\(\parcial{}{t} - \Delta\) v - \p'' |L|^2 |\nabla v|^2 +  \varphi \(\parcial{}{t} - \Delta\) |L|^2
\\ & \quad  - \frac1{\varphi} \<\nabla {\frak g}, \nabla \varphi\> + 2 \varphi |\nabla L|^2 + \frac3{2 \varphi} |L|^2 |\nabla\varphi|^2,
\end{align*}
where, exactly as in \cite{CaMi2}, we have used an inequality from \cite{EcHu} (combined with Kato's inequality
$| \nabla|L| | \le |\nabla L|$) to bound the last term in the first line.

By our hypotheses, we are working within a bounded domain of the ambient manifold $\oM$; in particular, all the curvatures of $\oM$ appearing  in the evolution formula (c) of Lemma \ref{Evol_gen} are bounded. Hence, arguing as in (6.12) of \cite{CaMi2}, we can find two positive constants $C_1$, $C_2$ so that
\begin{align}
\(\parcial{}{t} -  \Delta\) {\frak g}
\le & \ {\frak S} - |L|^2 \(\frac2{v \p'} + \frac{\p''}{\p'^2} - \frac{3}{2 \p}\) |\nabla \p|^2   - 2 \varphi \oH \vtr L^3     \nn \\
 &   
+ C_1 {\frak g} + C_2 \sqrt{\p {\frak g}} - \frac1{\varphi} \< \nabla {\frak g}, \nabla\varphi \>  , \lb{evolg1}
\end{align}
with 
$${\frak S} =- |L|^2 \p'\left( |L|^2 + \oRic(N,N) +  \frac{(n-1)}{f^2} \left[\frac{h'}{h}\right]'\right) v + 2 \varphi |L|^4,$$
where we have also used Lemma \ref{dvtn} and $|\nabla v| = |\nabla \p|/\p'$ to get, after rearranging and canceling terms, the inequality above.

Next, let us bound and/or rearrange the different terms in \eqref{evolg1}. First, from the definition of $\p$ in \eqref{fi}, it is easy to check
\be\lb{otrask}
\p' = \fracc{2 v}{(1- k v^2)^2} = 2 \frac{\p^2}{v^3}   \quad \text{ and } \quad \(\frac2{v \p'} + \frac{\p''}{\p'^2} - \frac{3}{2 \p}\)=  \frac{\varphi - v^2}{2 \varphi^2} .
\ee

Now we are in position to bound ${\frak S}$ as follows
\begin{align}
{\frak S} & = \(\frac2\varphi - \frac{\varphi'  v}{\varphi^2} \){\frak g}^2 - \( \oRic(N,N) +  \frac{(n-1)}{f^2} \left[\frac{h'}{h}\right]'\) \frac{\p'}{\p} v {\frak g} 
\le \(\frac2\varphi  - \frac{\varphi'  v}{\varphi^2}\){\frak g}^2 + K_0 {\frak g}. \lb{K1g}
\end{align}
Here $0<K_0 := {\frak R} \left\|\varphi' v/\varphi\right\|_{_\infty} \|f^{-2}\|_\infty  = 2 {\frak R} \left\|(1 - k v)^{-1}\right\|_{_\infty} \|f^{-2}\|_\infty \leq  {\frak R} \|f^{-2}\|_\infty$ where $\frak R$ is the constant coming from \eqref{defR}, and the last inequality is true by the choice of $k$.

Using $|\vtr L^3| \le |L|^3$ and Young's inequality  with $\varepsilon = k \varphi$ for $k$ as in \eqref{fi}:
\begin{align}
- 2 \varphi \oH \vtr L^3    \le 2 \varphi |\oH| |L|^3 = 2  |\oH| |L| {\frak g} \le \(k \varphi |L|^2 + \frac1{4 k \varphi} 4 \oH^2 \) {\frak g} = k {\frak g}^2 + \frac{\oH^2}{k \varphi} {\frak g}. \lb{boHtrL}
\end{align}

Plugging the expressions from \eqref{otrask}  to   \eqref{boHtrL} into \eqref{evolg1} and using Corollary \ref{p:uboH},   we reach the inequality in the statement 
for two positive constants $K_1$ and $K_2$.
\end{demo}

In \cite{CaMi2}, we managed to exploit a special symmetry of the problem in order to apply the maximum principle directly to the inequality corresponding to that from Lemma \ref{evolg} without having care of the boundary. However, our present setting lacks that symmetry, therefore, we need to have into account the effect of the boundary. To do so, we  consider another function $\psi(z) \frak g$, which gives us more {\it freedom} to get interior and boundary estimates of the heat operator acting on such a new function.

\section{Global curvature estimates}\lb{bL}

Before analyzing its behavior at the boundary (see Lemma \ref{Dztilg}), we have to deduce interior estimates for $\psi(z) \frak g$ (cf. Lemma \ref{evwg} below). A combination of the interior and boundary estimates for such an adhoc function will allow us to achieve global bounds for the curvature in Proposition \ref{boL}, which will close this section.

\begin{lema}\lb{evwg} Under the same hypotheses and notation than in Lemma \ref{evolg}, let us define $\wg = \psi(z) \  {\frak g}$, where $\psi$ is any real function satisfying
\begin{center}
$\psi$, $\psi'$ and $\psi''$ are bounded and $\psi>c$ for some constant $c > 0$.
\end{center}
Then there are positive constants $C_1$, $C_2$, $C_3$ such that 
\begin{align}
\(\parcial{}{t} -  \Delta\) {\wg} \le  - \< 2 \nabla \ln \psi + \nabla \ln \p, \nabla\wg\> -C_1 {\wg}^2 +
C_2 {\wg} + C_3  \sqrt{{\wg}}. \lb{ev_til_g}
\end{align}
\end{lema}

\begin{demo} As $\psi$ is positive, using part (b) of Lemma \ref{Evolrz} and Lemma \ref{evolg}, we have
\begin{align*}
\(\parcial{}{t} -  \Delta\) {\wg} & \leq  \wg  \,\frac{\psi'}{\psi} \left[\oH  \<N,\partial_z\> + \(\frac{u^2}{ f^2} - n\) \frac{f'}{f} \right] - \frac{\psi''}{\psi} \frac{u^2}{f^2} \wg  - 2 \< \nabla \psi, \nabla \fg\>  \\
& \quad 
+ \psi \(-k {\frak g}^2 +
K_1 {\frak g} + K_2 \sqrt{{\frak g}}  - \frac1{\varphi} \< \nabla{\frak g} , \nabla\varphi\>+\frac{v^2 - \varphi}{2 \varphi^3} |\nabla\p|^2 {\frak g}\).
\end{align*}

Now for the gradient terms we compute
\begin{align*}
- 2 \< \nabla \psi, \nabla \fg\> &=  \frac2{\psi}\<\nabla \psi, - \nabla \wg + {\frak g} \nabla \psi\>= - 2 \< \nabla \ln \psi, \nabla\wg\> +2 \frac{|\nabla\psi|^2}{\psi^2} \wg \\
-\frac{\psi}{\p}\< \nabla{\frak g} , \nabla\varphi\> &= \frac1{\p} \<\nabla\p, \fg \nabla\psi - \nabla \wg\> = \<\nabla \ln \p, \nabla \ln \psi\> \wg - \<\nabla \ln \p, \nabla\wg\>.
\end{align*}
Therefore
\begin{align*}
\(\parcial{}{t} -  \Delta\) \, \wg  & \leq   \left\{K_1 + \frac{\psi'}{\psi} \left[\oH \<N, \partial_z\> + \frac{f'}{f}\(\frac{1}{(f v)^2} -n \)\right] - \frac{\psi''}{\psi (f v)^2} + 2 \frac{|\nabla \psi|^2}{\psi^2}\right\} \wg  
\\ &   \quad -\frac{k}{ \psi} {\wg}^2 + K_2 \sqrt{\psi}\sqrt{\wg} 
- \<\nabla \wg, 2 \nabla \ln \psi + \nabla \ln \p\> +  \mc{D}\ \wg  
\end{align*}
with  $\mc{D} = \<\nabla \ln \psi, \nabla \ln \p\> + \ds \frac{v^2 - \p}{2 \p} \frac{|\nabla \p|^2}{\p^2}$.

\medskip

By explicit computation of $\nabla \p = \p' \nabla v$, it is easy to check that $\mc D$  is of order $|L|$. To compensate this, let us apply Young's inequality with $\varepsilon = k /\(2 \| f^2\|_\infty \)$:
\begin{align*}
\mc{D} &\le \frac{\| f^2\|_\infty}{2 k} \frac{|\nabla \psi|^2}{\psi^2} + \(\frac{k}{\| f^2\|_\infty} + \frac{v^2 - \p}{\p}\)\frac{|\nabla \p|^2}{2 \p^2}
\leq \frac{\| f^2\|_\infty}{2 k} \frac{\psi'^2}{\psi^2}  
\end{align*}
because $\nabla \psi = \psi'(z) \nabla z$,  \eqref{fi} implies $\fracc{k}{\| f^2\|_\infty} + \fracc{v^2}{\varphi} - 1 =  k \(\frac{1}{\| f^2\|_\infty} - v^2\)  \le 0$, and 
$|\nabla z| = |\frak t(z)| = |\<\frak t, \partial_z\>|   \leq 1.$

Plugging the above inequalities into the definition of $\wg$, we reach
\begin{align*}
\(\parcial{}{t} -  \Delta\) \, \wg  & \leq   \left[K_1 + \frac{|\psi'|}{\psi} \(|\oH|  + \frac{|f'|}{f}(n - 1)\) + \frac{|\psi''|}{\psi} + \frac{\psi'^2}{\psi^2}\(2 + \fracc{\| f^2\|_\infty}{2 k}\) \right] \wg  
\\ &  \quad -\frac{k}{ \psi} {\wg}^2 + K_2 \sqrt{\psi}\sqrt{\wg} 
- \<\nabla \wg, 2 \nabla \ln \psi + \nabla\ln \p \>. 
\end{align*}
Applying Corollary \ref{p:uboH} and using our hypotheses about $\psi$ and its derivatives, we deduce that all the coefficients of the different powers of $\tilde \fg$ in the above formula are bounded, which gives the  positive constants in the statement.
\end{demo}

\begin{lema} \lb{Dztilg}
Set $\psi(z) := [f(z)]^{-2m}$ for any $m >0$.  Under the same hypotheses that in the preceding section, 
at the boundary $\partial M_t$, for $t\in]0,T[$, one has 
$$\partial_z \wg =  2 \,\psi \,\p  \, \frac{f'}{f} \, \frak B,$$
 where 
$$ \frak B:= 3 (n-1) k_1 k_2 - (m + n - 2 + \p\, f^2) k_1^2 - (n -1)(m + 1 + \p\, f^2)k_2^2 - k_1 \oH.$$ \end{lema}

\begin{demo} First, let us compute 
\begin{align} \lb{DzL2}
\partial_z |L|^2 & = 2 k_1 \partial_z k_1 + 2 (n - 1) k_2 \partial_z  k_2 
 = 2 k_1 \partial_z H + 2(n -1) (k_2 - k_1) \partial_z k_2. 
\end{align}

To compute $\partial_z H$, recall that on the boundary we have $\frak t = \partial_z$ and $N = E_r$. Hence, along the boundary, at $t>0$
$$\oN_{\partial_t} E_r = \oN_{\partial_t} N = \nabla H = \frak t(H) \frak t = \partial_z H \partial_z,$$
where we have used formula (a) from Lemma \ref{Evol_gen}.  Using the flow equation, we also get
$$\oN_{\partial_t} E_r  = (\oH - H) \oN_N E_r = (\oH-H) \oN_{E_r} E_r =  (H - \oH) \frac{f'}{f} \partial_z,$$
which yields
\be \lb{DzH}
\partial_z H\big|_{\partial M} = (H -\oH) \frac{f'}{f}.
\ee

Next, taking into account \eqref{umbz}, \eqref{k2f2}, $\dot r|_{\partial M} = 0$ and that $\oN_{\partial_z} N|_{\partial M} = k_1 \partial_z$,
\begin{align} 
\partial_z k_2 \big|_{\partial M} &  = \partial_z\bigg(\frac{h'}{hf}\<N, E_r\> + \frac{f'}{f} \<N, \partial_z\>\bigg)\Big|_{\partial M} \nn
\\ & =  - \frac{h' f'}{h f^2} + \frac{h'}{hf} \<\oN_{\partial_z} N, E_r\>\big|_{\partial M} + \frac{f'}{f} \<\oN_{\partial_z} N, \partial_z\>\Big|_{\partial M} = \frac{f'}{f} (k_1 - k_2). \lb{Dzk2}
\end{align}

Now, substituting \eqref{DzH} and \eqref{Dzk2} in \eqref{DzL2}, we achieve
\begin{align} \lb{DzL22}
\partial_z |L|^2 &= 2 k_1 (H -\oH) \frac{f'}{f} - 2 (n - 1) \frac{f'}{f} \(k_2 - k_1\)^2\nn \\
&= 2 \frac{f'}{f} \(k_1^2 + (n -1) k_1 k_2 -k_1 \oH - (n -1)  k_2^2  + 2 (n -1) k_1 k_2 - (n -1) k_1^2\)\nn \\
&= 2 \frac{f'}{f} (- k_1 \oH - (n-2) k_1^2 + 3 (n-1) k_1 k_2 -  (n-1) k_2^2).
\end{align}

On the other hand, applying again \eqref{umbz} and $\oN_{\partial_z} N|_{\partial M} = k_1 \partial_z$, we obtain
\begin{align} \lb{Dzv}
\partial_z v &= \partial_z \(\<N , f E_r\>^{-1}\) = -v^2 \partial_z \(\<N , f E_r\>\)
=  - \frac1{f^2} \<N, E_r\> \partial_z f = - \frac{f'}{f^2}.
\end{align}

Then, from the definition of $\wg$ and substituting the explicit expression of $\psi$,
\begin{align*}
\partial_z \wg  &= \partial_z \(\psi(z) \p(v) |L|^2\) = \psi'  \fg + \psi\, \p'  |L|^2 \partial_z v + \psi \, \p \, \partial_z|L|^2 \\
& = -2m \frac{f'}{f} \psi \frak g - 2 \psi \frac{\p^2}{v^3} \frac{f'}{f^2} |L|^2 + \psi \, \p \, \partial_z|L|^2
\\ & \!\! = 2 \, \p\, \psi \frac{f'}{f}\Big(-(m + \p f^2) |L|^2  - k_1 \oH - (n-2) k_1^2 + 3 (n-1) k_1 k_2 -  (n-1) k_2^2\Big),
\end{align*}
where for the equality of the second line we have applied \eqref{otrask} and \eqref{Dzv}, and for the last equality we have used \eqref{DzL22}. 
Finally, substituting $|L|^2 = k_1^2 + (n - 1)k_2^2$ and rearranging terms, we reach the formula in the statement.
\end{demo}

\begin{prop}\label{boL}  Let $M_t$ be the solution of  \eqref{vpmf}  with initial condition $M$ in the setting $\frak{Eq}$ and satisfying the boundary condition \eqref{bocon}. If there are constants $\de$ and $\rho $ so that $\z >\de \ge r_t \geq \rho>0$, then we can find a positive constant $C_0 = C_0(V, \overline g,n,a,b,\rho,\de)$ such that $|L| \leq C_0$ on $M_t$.
\end{prop}

\begin{demo}
Observe that, when $\wg$ has the maximum within the interior, we can perform in the standard way a maximum principle argument for the inequality in Lemma \ref{evwg} (like in \cite{CaMi2}) to conclude that $\wg$ is bounded. Since $\psi$ and $\p$ are also bounded, we achieve the desired upper bound for $|L|$.

It remains consider the case of $\wg$ attaining the maximum  at the boundary, so that $\partial_z \wg|_{\partial M} \ge 0$. Notice that $\dot r|_{\partial M} =0$ which, by substitution in \eqref{k2f2}, gives 
$$|k_2| = \left|\fracc{h'(r)}{h(r) f(z)}\right| \le \fracc1{\min\{f(a), f(b)\}} \ds \left\| \fracc{h'}{h}\right\|_\infty=: \frak k_2.$$ 

Let us  assume  $|k_1| > \ell:= \max\{1,\frak k_2\}$. This allows us to estimate the quantity in Lemma \ref{Dztilg} as
\begin{align}
{\frak B} & <   3 (n-1) k_1 k_2 - k_1 \oH - m k_1^2 \le  3 (n-1) |k_1| \ell  + |k_1 \oH| - m k_1^2   \nn\\
& < |k_1|\, (3 (n-1) \ell + h_2 - m \ell),\nn
\end{align}
where $h_2$ is the constant coming from Corollary \ref{p:uboH}. If we choose $m \geq h_2 + 3(n -1)$, we obtain from Lemma \ref{Dztilg} that $\partial_z\wg |_{\partial M} < 0$, which contradicts the above assertion of $\partial_z \wg|_{\partial M} \ge 0$. 
 In conclusion, $|k_1| \le  \ell$, then $\wg = (k_1^2 + (n-1) k_2^2)\psi \p$ has an upper bound on $\partial M$ and thus $|L|$ is bounded. 
\end{demo}

Once we have uniform upper bounds for $|L|$, to get long time existence when the evolving manifold keeps away from the axis of rotation, we need that also all derivatives $|\nabla^k L|$ are bounded, which again will require a careful analysis of what happens on the boundary. We address this issue in the next section.

\section{The first singularities of the motion are produced at the axis of revolution}\lb{preLTE}

Here we prove the following result, which assures long time existence unless the evolving hypersurface reaches the axis of rotation.

\begin{teor}\lb{t:preLTE}
Let $M_t$ be the maximal solution of  \eqref{vpmf}, defined on  $[0,T[$, with initial condition $M$ in the setting ${\frak Eq}$ and satisfying \eqref{bocon}. Then either 
$$T=\infty \qquad \text{ or } \qquad\left\{ \begin{matrix} \ds \inf_{[0,T[} \min_{ M_t} r(\cdot, t) = 0  \\ \text{ or } \\ \z<\infty \text{ and }  \ds \sup_{[0,T[} \max_{ M_t} r(\cdot, t) = \z\end{matrix} \right. .$$ 
\end{teor}

\begin{demo}
First notice that, if the flow is defined on $[0,T[$, $\min_{M_t} r_t(x) > 0$ and $\max_{M_t} r_t(x) < \z$ for each $t\in[0,T[$.   Let us assume $\rho = \inf_{[0,T[} \min_{ M_t} r(\cdot, t) >0$ and $\de = \sup_{[0,T[} \max_{ M_t} r(\cdot, t) < \z $ ; then the goal is to show that the solution of the flow can be prolonged after $T$, which is a contradiction.

Since $r_t \in [\rho, \de]$, our evolving hypersurface remains within a bounded subset of the ambient manifold $\oM$; accordingly, we have bounds for the curvature $\oR$ and its covariant derivatives. Hence, following the same procedure of \cite{Hu84, Hu86}, we can find a constant $D_1 = D_1(n, \og, C_0, h_2)$ (where $h_2$ and $C_0$ are the constants coming from Corollary \ref{p:uboH} and Proposition \ref{boL}, respectively) such that 
\begin{align}
\parcial{}{t} |\nabla L |^2 &\leq \Delta |\nabla L |^2   +  D_1 (|\nabla L |^2 + 1). \lb{dtAm1}
\end{align}

 We define  
\be \lb{deff}
\Psi:= |\nabla L |^2 + \xi | L |^2,
\ee
for some positive constant $\xi$ to be specified later. 
 For the time derivative of $\Psi$, \eqref{dtAm1} yields
$$\parcial{\Psi}{t} \leq \Delta |\nabla L |^2  + D_1 \(|\nabla L |^2 + 1\) + \xi \parcial{}{t} |L |^2
\\  \leq \Delta \Psi + (D_1 - 2 \xi) |\nabla L|^2 + D_1 + \xi D_2,$$
where $D_2$ comes from bounding the curvature terms, $|\oH|$ and the different powers of $|L|$ in Lemma \ref{Evol_gen} (c).

If we choose $\xi \geq D_1$ and $D_3:= D_1 + \xi D_2$,   having into account \eqref{deff} and Proposition \ref{boL}, we deduce 
\begin{align*}
\parcial{\Psi}{t} & \leq \Delta \Psi -  \xi |\nabla L |^2 + D_3 \leq   \Delta \Psi  -  \xi \Psi + \xi^2 C_0^2 + D_3. 
\end{align*}
From here, a maximum principle argument ensures that $\Psi$  (then $|\nabla L|$) is bounded if it is bounded at the boundary.  To ensure that the latter indeed happens, we consider the equivalent flow equation \eqref{tafl}, take  derivatives with respect to $z$ and, after that, evaluate on $\partial M$ having in mind that   $\dot r = 0$, $\parcial{\dot{r}}{t} =0$, $|\dot{c}| = 1$ and $\partial_z |\dot{c}| = 2 \dot r f (\ddot r f + \dot r f') = 0$ on the boundary. Doing so, we obtain
\begin{align} \lb{8.4}
0 = \dddot{r} + (n+1) \frac{f'}{f} \ddot{r} + 2 (n-1) \frac{h'}{h}\frac{f'}{f^3} - \oH \frac{f'}{f^2}.
\end{align}
Again by Proposition \ref{boL},
$|L||_{\partial M}$ is  bounded which, combined with \eqref{k1f1}, gives a bound for $\ddot{r}|_{\partial M}$.  Thus \eqref{8.4} implies that  $\dddot r|_{\partial M}$ is bounded; this, again by \eqref{k1f1}, ensures that $|\nabla L|$ is bounded on the boundary, as we needed to show. Then  $|\nabla L|$ is bounded on the whole $M$,  and so is $\dddot{r}$. 

Next, we can substitute the solution $r(z,t)$ in \eqref{tafl} and see that it is a solution of the linear PDE
\begin{align}
\parcial{r}{t} &=   a(z,t) \ddot{r} + b(z,t), \lb{lpde}
\end{align}
where
\begin{align} \lb{ab}
a =  \frac{1}{1 + (\dot r f)^2}  \qquad 
b = \frac{f'}{f} \(\frac1{1 + (\dot r f)^2} + n\) - (n-1)  \frac{h'}{h f^2}+ \oH \frac{\sqrt{1 + (\dot r f)^2}}{f}.
\end{align}

Until now we have proved that  $r$, $\oH$, $\dot r$, $\ddot r$ and $\dddot r$ are bounded. Then from \eqref{tafl} it follows that also $\parcial{\dot r}{t}$ is bounded. Taking derivatives in \eqref{tafl} with respect to $z$ we obtain $\parcial{\dot r}{t}$ as a function depending on $z$, $r$, $\dot r$, $\ddot r$, $\dddot r$ and $\oH$, hence it is also bounded. Moreover, it follows from \eqref{oHt}, \eqref{I1oH}  and \eqref{muf} 
 that $\parcial{\oH}{t}$ is bounded if $r$, $\dot r$, $\parcial{r}{t}$ and $\parcial{\dot r}{t}$  are bounded (which we know is true) and $|M_t|$ is bounded from below by a positive constant. But this last condition follows from $|M_t| \ge \mathcal A_V$, the last being the $n$-volume of the hypersurface of minimum area enclosing a volume $V$  and with boundary orthogonal to and included between the hypersurfaces  $z=a$, $z=b$. Summing up, we conclude that $\parcial{b}{t}$ and $ \parcial{b}{z}$ are bounded.

Following the notation in \cite{rusas}  for the  H\"older norms, the bounds remarked before imply, by Theorem 5.4 page 322 in \cite{rusas},  that $|r|_{(3)}$ is bounded.  Repeating the argument (doing the standard bootstrapping argument), we have that, for every $m$, $|r|_{(m)}$ is bounded by some constant depending on $m$, then also $|\nabla^mL|$ is bounded, and arguing as in \cite{Hu84} we can continue the flow after $T$. 

\end{demo}

\section{Initial conditions giving  long time existence and convergence.} \label{volp}

\begin{teor}\lb{smallv} If, for an initial hypersurface $M$ in the setting ${\frak Eq}$, we have the following upper bound for the area $|M|$:
\be\lb{rAV}
|M| \leq  \frac{\min \left\{V, \vle(G)-V\right\}}{\| f^{-n}\|_{_\infty} \int_a^b f(z)^n dz},
\ee
  then the solution of \eqref{vpmf} satisfying \eqref{bocon}  is defined for all $t>0$, and there is a subsequence of times $t_n$ for which the corresponding solution converges to a revolution hypersurface of constant mean curvature in $\oM$.
\end{teor}

 In our setting $\frak Eq$, when $\z=\infty$, $\vle(G) =\infty$ and the hypothesis reduces to a uniform upper bound on the ratio $|M|/V$. When $\z<\infty$, $\vle(G)$ is finite, and $\vle(G)-V$ has the same right than $V$ to be called the volume enclosed by $M$, then the necessity to modify the hypothesis when $\vle(G)$ is finite is quite natural. 
\medskip 

\begin{demo} 
We can assume that our initial $M$ has non-constant $H$ (since, otherwise, it is a steady soliton of the flow \eqref{vpmf} and the statement follows trivially). Then Lemma \ref{Evol_gen} (b) implies
\be \lb{disp1}
|M_t| < |M| \quad \text{and } \quad r_2(t)<r_2 \quad \text{ for any } t>0,
\ee
where $r_2(t)$ is the upper bound of $r$ at time $t$ obtained by direct application of \eqref{defr2} using $|M_t|$ instead of $|M|$. 

Observe that, when $M_t$ is a graph (and we take $z=s$), the first inequality in \eqref{fVr1} becomes an equality which yields $r_{m}(t) \le r_1 \le r_{M}(t)$. Given $t_0>0$, we set 
$$\eps = |M|- |M_{t_0}| >0 \qquad \text{ and } \qquad \rho = \inf_{[t_0,T[}\{r_{m}(t)\}.$$
Now we use  the continuity of the function $F(\ell) = \omega_{n - 1} \min_{[a,b]} f(z)^n \int_{\ell}^{r_1} h(r)^{n-1} dr$  at $\ell = \rho$ to choose a $t_\eps \ge t_0$ so that $r_m(t_\eps)$ is close enough to $\rho$ in order to imply $F(r_m(t_\eps)) > F(\rho) - \frac{\eps}{2}$. Plugging the latter and $r_1 \leq r_M(t)$ into \eqref{vle>+} leads to
\begin{align}
|M| &= |M_{t_0}| + \eps \ge |M_{t_\eps}| + \eps > \omega_{n - 1} \min_{[a,b]} f(z)^n \int_{\rho}^{r_1}   h(r)^{n-1} dr + \frac{\eps}{2} \nn \\
&= \omega_{n - 1} \ \| f^{-n}\|_\infty^{-1}   \(\delta (r_1) - \int_0^{\rho}   h(r)^{n-1} dr\) + \frac{\eps}{2} \nn \\
&=     \frac{V }{\| f^{-n} \|_\infty \int_a^b f(z)^n dz} - \frac{\omega_{n - 1}}{\|f^{-n}\|_\infty}   \int_0^{\rho}   h(r)^{n-1} dr + \frac{\eps}{2}, \nn
\end{align}
where we have applied the definition of $r_1$ in \eqref{defR1}. Note that the above inequality is compatible with the hypothesis  \eqref{rAV}  only if $\rho>0$.
Hence 
\be\lb{lgbr}
r_t\ge \rho >0 \quad \text{ for every } \quad  t\in[t_0,T[. 
\ee

On the other hand, the quantity  $\vle(G)$ can be written (cf. \eqref{fVr1}) as
$$\vle (G) = \omega_{n-1}\(\int_a^b f(z)^n dz \) \(\int_0^\z h(r)^{n-1} dr\),$$
which gives
$$\delta(\z) = \frac{1}{\omega_{n-1}} \frac{\vle(G)}{\int_a^b f(z)^n dz}.$$
Next, \eqref{defr2} together with the hypothesis \eqref{rAV} imply
\begin{align*}
\delta(r_2)  & = \frac1{\omega_{n-1}}  \(|M| \|f^{-n}\|_\infty + \frac{V}{\int_a^b f(z)^{n}dz}\)   \le 
\frac1{\omega_{n-1}} \( \frac{\vle(G)-V}{\int_a^b f(z)^n dz}+ \frac{V}{\int_a^b f(z)^{n}dz}\) \nn \\
&
= \delta(\z).  
\end{align*}
Accordingly, $r_2 \le \z$ and, thanks to \eqref{disp1}, we have $\eta=r_2- r_2(t_0) > 0$. Using Proposition \ref{borM} for any fixed time $t\ge t_0$ and the definition of $r_2(t)$ combined with the decreasing of area under the flow, we reach
\be\lb{ugbr}
r_t < r_2(t) \le r_2(t_0) = r_2 -\eta =:\de \le \z -\eta <\z.
\ee

Then, from \eqref{lgbr} and \eqref{ugbr} we conclude, because of  Theorem \ref{t:preLTE}, that the solution of \eqref{vpmf} is defined on $[0,\infty[$; hence $r$ is bounded uniformly from above and below on $[0,\infty[$. After the results of section \ref{graph}, it is clear that  $\dot{r}$ remains bounded all the time. 
In addition, the proof of Theorem \ref{t:preLTE} shows that  $|\nabla^jL|^2$ is uniformly bounded for every $j\ge 0$.  Once we have all these bounds, it follows from \eqref{k1f1} and \eqref{k2f2} (taking $z=s$) that all the derivatives of $r$ are bounded on $[0,\infty[$.  

We are now in position to  apply Arzel\`a-Ascoli Theorem to ensure the existence of a sequence of maps  $ r_{t_i}$ satisfying \eqref{tafl} which $C^\infty$-converges to  a smooth  map  $r_\infty : [a,b]\flecha \re^+$ also solving \eqref{tafl}. A standard argument like in \cite{CaMi1} proves that the limiting hypersurface $ M_\infty =(z, r_\infty(z),u)$
has constant mean curvature.
\end{demo}

Here we give a final remark for those readers familiar to \cite{CaMi2},
who may wonder why the above proof is not as short as that of the
corresponding result (namley, Theorem 12) in \cite{CaMi2}. The reason is
that there is a \lq\lq typo'' in the hypothesis on the inequality
satisfied by $|M|$: where it says \lq\lq $\le$'', it should say \lq\lq
$<$ ''. To attain the same result using the weaker assumption \lq\lq
$\le$'' as in the theorem above, one needs to obtain finer estimates
like in the previous proof. On the other hand, the proof of the
convergence of the sequence in \cite[Theorem 12]{CaMi2} has an issue, in
fact, what is actually proved there is the existence of a convergent
subsequence.

\begin{appendix}

\section{Appendix - Proof of Lemma \ref{Evolrz}}

First, we shall obtain the evolution of $r$. 
\begin{align} 
\parcial{r}{t} = (\oH-H) \< \oN r,N\> \us{\eqref{gradr}}{=} (\oH-H) \fracc{u}{f^2}. \lb{drdt}
\end{align}
Since $\nabla_{\frak t} \frak t = 0$ and $E_i(r) = 0$ 
 (as can be easily deduced from Lemma \ref{Conn}), we have
\begin{align}
\Delta r = \frak t \frak t (r) - \sum_{i=2}^{n}  \nabla_{E_i}  E_i (r). \lb{Delr}
\end{align}

We compute, using \eqref{gradr} and \eqref{t},
\be
\frak t(r) = \<\oN r, \frak t\> = \frac1{f} \<E_r, \frak t\>  =  -\frac1{f} \<N, \partial_z\> \lb{tr}
\ee
and, using \eqref{t} and Lemma \ref{Conn} repeatedly,
\begin{align}
\frak t \frak t (r) &= \frac{f'}{f^2} \<\frak t, \partial_z\> \<N, \partial_z\> - \frac1{f}\<\oN_{\frak t} N, \partial_z\> - \frac1{f} \<N, \oN_{\frak t} \partial_z\> \nn
\\ & = \frac{f'}{f^2} \<N, E_r\> \<N, \partial_z\> - \frac1{f} \<t, \partial_z\> k_1 - \frac1{f} \<t, E_r\> \<N, \oN_{E_r} \partial_z\>\nn
\\ & =  \frac{f'}{f^3} u \<N, \partial_z\> -  \frac1{f} \<N, E_r\> k_1 + \frac1{f} \<N, \partial_z\> \<N, \frac{f'}{f} E_r\>\nn
\\ & =   2 \frac{f'}{f^3} u \<N, \partial_z\> - \frac{u}{f^2} k_1. \lb{ttr}
\end{align}

On the other hand, using \eqref{Dij}, \eqref{t} and \eqref{k2f2},
\begin{align}
(\nabla_{E_i} E_i)(r)  & = \<\oN_{E_i} E_i, \frak t\> \frak t (r) = \(\frac1{f} \frac{h'}{h} \<\frak t, E_r\> + \frac{f'}{f} \<\frak t, \partial_z\>\)\frac1{f} \<N, \partial_z\> \nn
\\ & \!\! =  - \frac{h'}{f^2 h} \<N, \partial_z\>^2 + \frac{f'}{f^2} \<N, E_r\> \<N, \partial_z\>
= \frac1{f^2} \(u \, k_2 - \frac{h'}{h}\). \lb{Ei2r}
\end{align}
Now, substituting \eqref{ttr} and \eqref{Ei2r} in \eqref{Delr}, we have
\be
\Delta r = 2 \frac{f'}{f^3} u \<N, \partial_z\> + (n -1) \frac{h'}{f^2 h} - \frac1{f^2} H u. \lb{Delr2}
\ee 
Therefore, for any function $\phi: \re \fle \re$, we have
\begin{align}
\(\parcial{}{t}-\Delta\) &\phi(r) = \phi'\(\parcial{}{t} - \Delta \)r - \phi'' |\nabla r|^2 \nn \\
&= \phi'\left( \oH \frac{u}{f^2} - 2 \frac{f'}{f^3} u \<N, \partial_z\> - (n - 1) \frac{h'}{f^2 h} \right) +  \frac{\phi''}{f^2} \(\frac{u^2}{f^2} - 1\),
\end{align}
where the second equality follows plugging \eqref{Delr2} into \eqref{drdt}, and we have also used that 
$|\nabla r|^2 = |\frak t(r)|^2  =  \fracc1{f^2} \<N, \partial_z\>^2 .$

With the goal to prove part (b) of Lemma \ref{Evolrz}, we first need the evolution equation of the axial coordinate $z$.
\begin{align}
\parcial{z}{t} &= (\oH-H) \<\oN z, N\> = (\oH - H) \<N,\partial_z\>. \lb{dzt1}
\end{align}

From Lemma \ref{Conn} it is easy to compute that
$$\Delta z = \frak t \frak t(z)  + \sum_{i=2}^n  \nabla_{E_i} E_i (z),$$
\begin{align*}
\text{with } \qquad \frak t (\frak t z) &= \frak t \(\<\frak t, \partial_z\>\) = \<\oN_{\frak t} \frak t, \partial_z\> + \<\frak t, \oN_{\frak t} \partial_z\>  
\\ & =  - k_1 \<N, \partial_z\> + \<\frak t, E_r\> \big<\frak t, \oN_{E_r} \partial_z \big>= - k_1\<N, \partial_z\> + \frac{f'}{f} \<N, \partial_z\>^2,
\end{align*}
and, using again \eqref{Dij}, \eqref{t} and \eqref{k2f2},
\begin{align*}
(\nabla_{E_i} E_i)(z) & = \<\oN_{E_i} E_i, \frak t\> \frak t(z) =  -\(\frac{f'}{f} \<\frak t, \partial_z\>^2 + \frac{h'}{f^2 h} \<\frak t, \partial_r\> \<\frak t, \partial_z\>\)
\\ & = -\frac{f'}{f}\<N, E_r\>^2 + \frac{h'}{hf} \<N, E_r\> \<N, \partial_z\>
\\ &  =  - \frac{f'}{f} \<N, E_r\>^2 + \<N, \partial_z\> k_2 - \<N, \partial_z\>^2 \frac{f'}{f} = - \frac{f'}{f} +
\<N, \partial_z\> k_2 .  
\end{align*}

The above computations lead to
\be
\Delta z = - H \<N, \partial_z\> + \frac{f'}{f} \(n-\frac{u^2}{f^2}\)  \lb{Naz}
\ee

Finally, for any $\psi : \re \fle \re$, we get from \eqref{dzt1} and \eqref{Naz}
\begin{align*}
\(\parcial{}{t} - \Delta\) \psi(z) &= \psi' \(\parcial{}{t} - \Delta\) z - \psi'' |\nabla z|^2
\\ & = \psi \left(\oH \<N, \partial_z\> + \frac{f'}{f}\(\frac{u^2}{f^2} - n\)\right) - \psi''\frac{u^2}{f^2}. 
\end{align*}

\section{Appendix  - A hypersurface in the setting $\frak Eq$ with negative averaged mean curvature.}

For the cases with $\z <\infty$ there is a simple argument  showing that, in the setting $\frak Eq$, there must be revolution hypersurfaces inside $G$ with boundary orthogonal to $\partial G$ and $\oH <0$. In fact, let us suppose that $\oH>0$. The metric of $\oM$ can also be written, taking $(\z,u)$ instead of $(0,u)$ as the center of the spherically symmetric $(\mc S,\sigma)$,  as $dz^2 + f(z)^2 d\bar r^2 + f(z)^2 h(\z-\bar r)^2 g_{\Sn}$. With the metric written this way, the old axis $\mc A_-$ will be called $\widetilde{\mc A}_+$ now, and $M$ will be given as generated by the graph over $\widetilde{\mc A}_+$ of the function $\widetilde{r}(z) = r(\z-z)$. Now the domain bounded by $M$ will be $G-\Omega$,  and, as we consider positive the orientation given by the normal pointing outward, the positive orientation is now reversed respect to the original one, which gives  $\oH <0$.

Anyway, we give here an explicit case with $\oH<0$ for $\z<\infty$ and another for $\z=\infty$.

First, we consider the case (C2) of a revolution hypersurface inside a spherical crown in $\re^3$. As axis of symmetry we choose the axis $y$, and as generating curve we take a part of the cycloid 
$$(x(s), y(s)) = \(2 s -  \sin(s/2) + 2 \pi,  2 - \cos(s/2)\),$$
 which, written with the coordinates $(z,r)$ used to describe $(\oM,\og)$ (see the picture on the left), is 
$$(z(s), r(s)) =\(\sqrt{x(s)^2 + y(s)^2}, \arctan \(x(s)/y(s)\)\).$$

More precisely, we pick the portion $(z(s), r(s))$ for $s\in [s_1, s_2]$, where $s_1=4.33453$ and $s_2=12.7571$ are two consecutive values of $s$ satisfying $\dot r(s_i)=0$ (which guarantees that the revolution hypersurface generated by this curve is orthogonal to the boundary of the spherical crown $G$ between the spheres $z=z(s_1)$ and $z=z(s_2)$). If now apply the formula \eqref{oHt}  for $\oH$ , we obtain $\oH =  - \fracc{\omega_{n-1}}{|M|} 1.55553 <0$.
\includegraphics[scale=0.42]{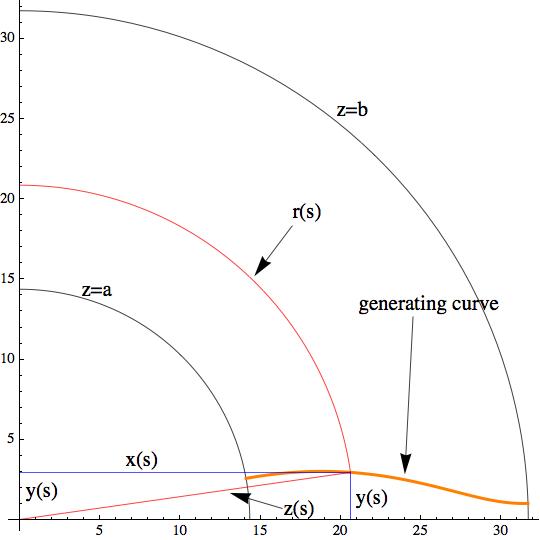}
\includegraphics[scale=0.5]{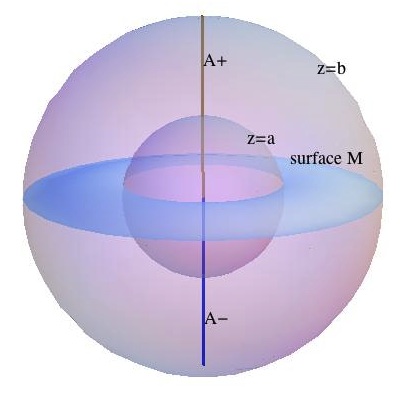}

Although explicit, the fact that $\z<\infty$ and the remark done at the beginning of this  appendix make this example not too much interesting. However we shall take the same expressions of $z(s)$, $r(s)$, $s_1$ and $s_2$ to obtain an example in the case (C5) (which obviously corresponds to $\z=\infty$) of a revolution hypersurface between two parallel horospheres in the hyperbolic space of dimension $3$. Using again formula \eqref{oHt} we obtain $\oH =  - \fracc{\omega_{n-1}}{|M|} \  9.72488 \ 10^{24} <0$.

\end{appendix}

\medskip

\noindent {\bf Acknowledgments.}	This research was partially supported by grants  DGI(Spain) and FEDER Project MTM2007-65852 and the Generalitat Valenciana Project Prometeo 2009/099. The first author was partly supported by the Leverhulme Trust.

{\footnotesize

\bibliographystyle{alpha}

}




\end{document}